\input amstex
\input epsf

\TagsOnRight
\magnification=\magstephalf

\def\E{{\,{\Bbb E}\,}}

\def\P{{\,{\Bbb P}\,}}
\def\di={\overset{\text{${\Cal D }$}}\to =}
\define\llapp#1{\hbox{\llap{${#1}$}}}

\define\rlapp#1{\hbox{\rlap{${#1}$}}}

\define\na{\,\, {\raise.4pt\hbox{$\shortmid$}}{\hskip-2.0pt\to}\, \, }

\define\<{\bigl <}\define\8{\left<}\define\9{\right>}
\define\>{\bigr>}
\define\ovln#1{\,{\overline{\!#1}}}

\define\PM{^{\2\prime}}

\redefine\D #1{{\<D{#1},{#1}\>}}
\define\norm#1{\left\|#1\right\|}
\define\nnorm#1{\bgl\|#1\bgr\|}
\define\snorm#1{\bigl\|\2#1\2\bigr\|}
\define\ssqrt{^{1/2}}


\define\weht#1{\bigl(\2\e+\nnorm{#1}\t\2\bigr)}
 \define\leei#1#2{\log\E e^{i\8\right.\!{#1},{#2}\!\left.\9}}
 \define\eei#1#2{\E e^{i\8\right.\!{#1},{#2}\!\left.\9}}
 \define\lee#1#2{\log\E e^{\8\right.\!{#1},{#2}\!\left.\9}}
 \define\ee#1#2{\E e^{\8\right.\!{#1},{#2}\!\left.\9}}
\define\xxxi#1{{\overline\xi}_{#1}(h)}

\define\oxh#1#2{{\overline\xi_{#1}}(h_{#2})}
\define\ooxh#1#2{{\overline\xi{}_{#1}\PM}(h'_{#2})}
\define\xxi#1{\xi_#1}
\define\4{\hskip1pt}

\define\wthj#1{\widetilde h_{#1}}

\define\twhtj#1{{\th\4\nnorm{\wthj#1}\tau}}

\define\xwthj#1{{\overline\xi_{{#1}}}(\wthj{#1})}
\define\xhpj#1#2{{\overline\xi{}\PM_{{#1}}}(h'_{#2})}

\define\detDj#1{(2\4\pi)^{d/2}\,(\det D_{#1})\ssqrt}

\define\2{\hskip.5pt}
\define\bgr#1{\4\bigr#1}
\define\bgl#1{\bigl#1\4}
\define\3{\hskip-10pt}
\define\R#1#2{\Cal R_d\bgl(#1,#2\bgr)}
\define\RB#1#2{\Cal R_d\Bgl(#1,#2\Bgr)}
\define\RBM#1#2{\Cal R_{d-1}\Bgl(#1,#2\Bgr)}

\define\pfrac#1#2{{}^{#1}\hskip-.5pt\!\big/\!\hskip.5pt{}_{#2}}
\define\sfrac#1#2{{#1}\big/{#2}}
\define\ffrac#1#2{\raise.5pt\hbox{\eightpoint$\4\dfrac{\,#1\,}{\,#2\,}\4$}}
\define\Bgr#1{{\eightpoint\4\raise.5pt\hbox{{$\Bigr#1$}}}}
\define\Bgl#1{{\eightpoint\raise.5pt\hbox{{$\Bigl#1$}}\4}}
\define\step#1{\raise9pt\hbox{\eightpoint$\tsize{#1}$}}

\documentstyle {amsppt}
\hoffset  0.25  true cm
\voffset  1.5 true cm

\hsize=16 true cm
\vsize=22 true cm
\loadmsbm
\def\t{\theta}

\def\D{\Delta }
\def\ve{\varepsilon }

\def\vk{\varkappa}
\def\ga{\gamma }

\def\l{\lambda }

\def\s{\sigma }

\def\s{\sigma}
\def\t{\theta}

\def\={\overset{ \text{def} }\to =}

\def\R{\Bbb R}

\def\Z{\Bbb Z}

\def\I{\Bbb I}

\define\e{\operatorname{e}}

\def\fs#1#2#3{ {#1}_{#2},\dots ,  {#1}_{#3}    }
\def\fsu#1#2#3{ {#1}_{#2}+\dots + {#1}_{#3}    }

\define\bbinom#1#2{\raise.5pt\hbox{\eightpoint$\dsize\4
\binom{\raise-1.3pt\hbox{${\ssize\,#1\,}$}}
{\raise1.3pt\hbox{${\ssize \,#2\,}$}}\4$}}

\documentstyle{amsppt}

\topmatter
\title Bounds for tail probabilities of martingales
using skewness and kurtosis
 \endtitle
\rightheadtext{ Skewness and curtosis  } \leftheadtext{V. Bentkus}
\author V. Bentkus\footnotemark"$^{1}$" and  T. {Ju\v skevi\v cius}$^{1}$
\endauthor

\footnotetext"$^{1}$"{ Vilnius Institute of Mathematics and
Informatics. Supported by grant ???}

\address
Vilnius institute of mathematics and informatics,\ 
Akademijos 4,\ 
LT-08663  Vilnius
\endaddress

\comment
and \smallskip
Vilnius pedagogical university \smallskip
Studentu 39\smallskip
LT-08106  Vilnius \smallskip
Lithuania  \smallskip
\endcomment

\email
bentkus\@ktl.mii.lt, ?????????????
\endemail

\date January 2008
\enddate
\keywords
 skewness, kurtosis,
 Hoeffding's
inequalities,
 sums of independent random variables, martingales,
bounds for tail probabilities, Bernoulli random variables, binomial
tails, probabilities of large deviations, method of bounded
differences
  \endkeywords
\subjclass
 60E15
\endsubjclass

 \abstract
Let $M_n= \fsu X1n$ be a sum of independent random variables such that
$ X_k\leq 1$,  $\E X_k =0$ and $\E X_k^2=\s_k^2$ for all $k$.
Hoeffding 1963, Theorem 3, proved that
$$\P\{M_n \geq nt \}\leq H^n(t,p ) ,\quad H(t,p)= \bgl( 1+qt/p\bgr)^{p +qt}
\bgl( {1-t}\bgr)^{q -qt}$$
with
$$q=\ffrac 1{1+\s^2},\quad p=1-q, \quad \s^2 =\ffrac { \s_1^2+\dots
+\s_n^2}n,\quad 0<t<1.$$
Bentkus 2004 improved Hoeffding's
inequalities using binomial tails as upper bounds.
Let $\ga_k =\E X_k^3/\s_k^3$ and $ \vk_k= \E X_k^4/\s_k^4$ stand for the skewness
and kurtosis of $X_k$.
In this paper we prove (improved) counterparts of the Hoeffding inequality
replacing $\s^2$ by certain functions of $\fs \ga 1n$ respectively $\fs \vk 1n$.
Our bounds extend to a general setting where $X_k$ are martingale
differences, and they
 can combine
the knowledge of skewness and/or kurtosis and/or variances of ~$X_k$.
 Up to  factors
bounded by  $e^2/2$ the bounds are final.
 All our results are new since no inequalities incorporating
skewness or kurtosis control so far are known.
 \endabstract

\endtopmatter

\document


\head  1. Introduction and results
\endhead

In a celebrated paper of Hoeffding 1963 several inequalities for sums of bounded random
variables were established. For improvements of the Hoeffding inequalities and
related results see, for example,   Talagrand ~1995, McDiarmid 1989,
Godbole and Hitczenko~ 1998,  Pinelis ~1998--2007,
Laib ~1999, B ~2001--2007,
van de Geer ~2002, Perron~ 2003,
BGZ ~2006--2006, BGPZ ~2006,
BKZ ~2006, ~2007, BZ 2003, etc.
Up to certain constant factors, these improvements are close to the final
optimal inequalities, see ~B ~2004, ~BKZ ~2006.
 However so far no  bounds taking into account  information related to
 skewness and/or kurtosis are known, not to mention
 certain  results related to symmetric random variables, see  BGZ 2006, BGPZ 2006.
In this paper we prove general and optimal counterparts of
Hoeffding's 1963 Theorem 3, using assumptions related to skewness and/or  kurtosis.

\bigskip

Let us recall Hoeffding's 1963 Theorem 3.
Let $M_n =\fsu X1n$ be a sum of independent random variables such that
 $ X_k\leq 1$, $\E X_k=0$, and $\E X_k^2=\s_k^2$ for all $k$.
Write
$$ \s^2=\ffrac {\s_1^2+\dots +\s_n^2}n, \quad p=\ffrac {\s^2} {1+\s^2},\quad q=1-p .$$
Hoeffding ~1963, Theorem 3, established the inequality
$$\P\{M_n \geq nt \}\leq H^n(t,p ) ,\quad H(t,p)= \bgl( 1+qt/p\bgr)^{p +qt}
\bgl( {1-t}\bgr)^{q -qt}  \tag 1.1$$
assuming that $ 0<t <1$.
One can rewrite $ H^n(x,p ) $ as
$$ H^n(x,p )=\inf_{h>0} \exp\{ -h nt \}\E \exp\{hT_n\},$$
where $T_n=\fsu \ve 1n $ is a sum of $n$
independent copies of a  Bernoulli random variable, say ~${\ve=\ve(\s^2)}$,
such that
 $$\P\{\ve =-\s^2 \}=q,\quad   \P\{\ve=1\}=p,\quad \E \ve^2 =\s^2.\tag 1.2 $$
Using the shorthand $x=nt$, we can rewrite the Hoeffding result as
$$\P\{M_n \geq x \}\leq \inf_{h>0} \e^{ -h x }\E \e^{hT_n}, \tag 1.3$$
In B 2004 the inequality (1.3) is improved to
$$\P\{M_n \geq x \}\leq \inf_{h<x} (x-h)^{-2} \E (T_n-h)_+^2
.\tag 1.4 $$
Actually, inequalities  (1.1), (1.3) and (1.4) extend to  cases where $M_n$ is a martingale or even super-martingale,
see B~ 2004 for a proof.
In  the case of (1.1) and (1.3) this was noted already by Hoeffding 1963.

 The right hand side of (1.4) satisfies
$$ \inf_{h<x} (t-h)^{-2} \E (T_n-h)_+^{2}\leq \ffrac {\e^2}2\P \{ T_n\geq x \}, \quad \e=2.718... \tag 1.5 $$
for integer $x\in \Z$. For non-integer $x$ one has to interpolate the probability
$\log$-linearly, see ~B ~2004 for details. The right-hand side of (1.4) can be given explicitly
as a function of $x$, $p$ and $n$,
see~ BKZ 2006, as well as Section 2 of the present paper. To have
bounds as tight as possible is essential for statistical applications, like those in audit,
see BZ 2003.

Our intention in this paper
 is to develop methods leading to counterparts of (1.1), (1.3) and ~(1.4) such that
 information related to the skewness and kurtosis
 $$\ga_k= \ffrac {\E ( X_k-\E X_k)^3}{\s_k^3},\quad \quad \vk_k= \ffrac {\E (X_k-\E X_k)^4}{\s_k^4}\tag 1.6$$
 of $X_k$ is taken into account
 (in this paper we define  $\ga_k=\infty$ and $\vk_k=1$
 if $\s_k=0$).
  All our results hold in general martingale setting.

All known proofs of inequalities of type (1.3) and (1.4) start with an application of
Chebyshev's inequality. For example, in the case of (1.4)
we can estimate
$$\P\{M_n \geq x \}\leq \inf_{h<x} (x-h)^{-2} \E (M_n-h)_+^2
\tag 1.7 $$
since the indicator function $t\mapsto \I \{ t\geq x\}$ obviously satisfies
$\I \{ t\geq x\}\leq (x-h)^{-2}  (t-h)_+^2$ for all ~${t\in \R}$.
The further proof of (1.4) consists in showing that
$\E (M_n-h)_+^2\leq \E (T_n-h)_+^2$ for all ~${h\in \R}$. We would like to emphasize that
all our proofs are optimal in the sense that no further improvements are possible
in estimation of $\E (M_n-h)_+^2$. Indeed, in the special case $M_n=T_n$
the inequality $\E (M_n-h)_+^2\leq \E (T_n-h)_+^2$ turns into the equality
$\E (T_n-h)_+^2= \E (T_n-h)_+^2$.

 In view of (1.7) it is  natural to introduce and to study
 transforms $G\mapsto G_\beta$
 of survival functions $G(x)=\P \{X\geq x\}$
 of the type
 $$G_\beta (x) = \inf_{h<x} (x-h)^{-\beta}\E (X-h)_+^\beta ,  \quad  \beta>0,\tag 1.8$$
defining $G_0=G$ in the case $\beta =0$. See Pinelis ~1988, 1989, B~ 2004, ~BKZ 2006 for related known results.

The paper is organized as follows. In the Introduction we provide necessary definitions and formulations
of  our results, including their versions for sums of martingale differences.
In Section 2 we recall a description of the transform $G\mapsto G_2$ of binomial survival
functions---our bounds are given using $G_2$.
Section 3 contains proofs of the results.

\bigskip
Henceforth $M_n=\fsu X1n$ stands for  a martingale sequence such that the differences ~$X_k$
are uniformly  bounded (we set $M_0=X_0=0$).
 Without loss of generality we can assume that the bounding constant is $1$,
that is, that  $X_k\leq 1$. Let $\Cal F_0 \subset \Cal F_1\subset \dots \subset \Cal F_n$
be a related sequence of $\sigma$-algebras
such that $M_k$ are $\Cal F_k$-measurable.
Introduce the conditional variance $s_k^2$, skewness ~$g_k$ and kurtosis $c_k$ of $X_k$ by
$$s_k^2= \E \bgl(X_k^2\, \Bgl|\, \Cal F_{k-1}\bgr),\quad
g_k= \E \bgl(  {X_k^3}\, \Bgl|\, \Cal F_{k-1}\bgr)/{s_k^3} ,\quad
c_k= \E \bgl( {X_k^4}\, \Bgl|\, \Cal F_{k-1}\bgr)/{s_k^4}.\tag 1.9
$$
Note that $s_k^2$, $g_k$, $c_k$ are $\Cal F_{k-1}$-measurable random variables.

\bigskip

{\bf Remark} {\bf 1.1}. We prove our results using (1.4) for martingales. It is proved in B ~2004 that
all three inequalities ~(1.1), ~(1.3) and
~(1.4) hold with $\s^2=(\s_1^2+\dots +\s_n^2)/n$
if $M_n$ is a martingale with differences ~${X_k\leq 1 }$ such that the conditional variances $s_k^2$
satisfy $s_k^2\leq \s_k^2$ for all $k$.

\bigskip

It is easy to check that Bernoulli random variables $\ve=\ve(\s^2)$ of type (1.2) have  variance $\s^2$ and
skewness ~$\ga$
related as
$$\ga = \s -\ffrac 1\s ,\quad \s^2 = u^2(\ga ),\quad \text{where} \ u(x) =\sqrt{1+\ffrac {x^2}4}-\ffrac x2
 .\tag 1.10 $$

\bigskip

\proclaim{Theorem 1.2} Assume that the differences $X_k$ of a martingale $M_n$
satisfy ~${X_k\leq 1}$, and that the conditional skewness $g_k$
of $X_k$ are bounded from below by some non-random $\ga_k$, that is, that
$$ g_k \geq \ga_k,\quad k=1,2,\dots , n. \tag 1.11$$
Then $(1.3)$ and $(1.4)$ hold
with $T_n$ being  a sum of $n$ independent copies of a Bernoulli random variable ~${\ve=\ve(\s^2)}$
of type $(1.2)$ with
skewness $\ga$ and variance $\s^2$ defined by
$$\ga= \ffrac 1{\sqrt{n}}\ffrac{\ga_1 u(\ga_1) +\dots + \ga_n u (\ga_n) }{\sqrt{u^2(\ga_1)+ \dots +u^2(\ga_n)}},
\quad \s^2 = \ffrac{ u^2(\ga_1)+\dots +u^2(\ga_n) }n. \tag 1.12 $$
In the special case where all $\ga_k$ are equal, $\ga_1=\dots =\ga_n=\ga$, the Bernoulli random variable has
skewness $\ga$ and variance $\s^2= u^2(\ga)$.
\endproclaim

\bigskip

 It is easy to see that Bernoulli random variables $\ve=\ve(\s^2)$ of type (1.2) have  variance $\s^2$ and
kurtosis ~$\vk$
related as
$$
 \vk  = \ffrac 1{\s^2}-1+\s^2 ,\quad 2\s^2=  {\vk+1\pm \sqrt {(\vk +1)^2-4}}
 .\tag 1.13 $$
In particular
$$
 \s^2\leq v(\vk ),
\quad \text{where} \ 2v(t) =t+1+ \sqrt {(t +1)^2-4}.\tag 1.14$$

\bigskip
\proclaim{Theorem 1.3}
 Assume that the differences $X_k$ of a martingale $M_n$
satisfy ~${X_k\leq 1}$, and that the conditional kurtosis $c_k$
of $X_k$ are bounded from above by some non-random $\vk_k$, that is, that
$$ c_k \leq \vk_k,\quad k=1,2,\dots , n. \tag 1.15$$
Then $(1.3)$ and $(1.4)$ hold
with $T_n$ being  a sum of $n$ independent copies of a Bernoulli random variable ~${\ve=\ve(\s^2)}$
of type $(1.2)$ with
kurtosis  $\vk$ and variance $\s^2$ defined by
$$\vk  = \ffrac 1{\s^2}-1+\s^2,
\quad \s^2 = \ffrac{ v(\vk_1)+\dots +v(\vk_n) }n,\tag 1.16 $$
where the function $v$ is given in $(1.14)$.
In the special case where  $\vk_1=\dots =\vk_n=\vk$, the Bernoulli random variable has
kurtosis  $\vk$ and variance $\s^2= v(\vk)$.
\endproclaim

\bigskip
The next Theorem 1.4 allows to combine our knowledge about variances, skewness and kurtosis. Theorems 1.2, 1.3
and (1.4), (1.4) for martingales (see Remark 1.1)
 are special cases of Theorem 1.4 setting in various combinations $\s_k^2=\infty$, $\ga_k=-\infty$, $\vk_k=\infty$.

\bigskip
\proclaim{Theorem 1.4}
 Assume that the differences $X_k$ of a martingale $M_n$
satisfy ~${X_k\leq 1}$, and that their
 conditional variances $s_k^2$, skewness $g_k $ and
 kurtosis $c_k$ satisfy
$$s_k^2\leq \s_k^2,\quad g_k \geq \ga_k,\quad
 c_k \leq \vk_k,\quad k=1,2,\dots , n \tag 1.17$$
with some non-random $s_k^2\geq 0$, $g_k\geq -\infty$ and $1\leq c_k\leq \infty$.
Assume that  numbers $\alpha_k^2$ satisfy
$$\alpha_k^2 \geq \min\{ \s_k^2,\,  u^2(\ga_k),\,  v(\vk_k)\}.$$
Then $(1.3)$ and $(1.4)$ hold
with ~${T_n}$ being  a sum of $n$ independent copies of a Bernoulli random variable ~${\ve=\ve(\s^2)}$
of type $(1.2)$ with
$$ \s^2 = \ffrac{ \alpha_1^2 +\dots +\alpha_n^2 }n,$$
where  functions $u$ and  $v$ are defined  in $(1.10)$ and $(1.14)$ respectively.
\endproclaim

\bigskip
{\bf Remark } {\bf 1.5}. All our  inequalities can be extended to the case where $M_n$ is a super-martingale.
Furthermore, their maximal versions hold, that is, in the left hand sides of these inequalities we
can replace $\P\{M_n\geq x\}$ by $\P\left\{\max\limits_{1\leq k \leq n}M_k\geq x\right\}$.

\bigskip
{\bf Remark } {\bf 1.6}. One can estimate the right hand sides of our inequalities
using Poisson distributions. In the case of Hoeffding's functions this is done by Hoeffding ~1963.
In notation of (1.1) his bound is
$$ H^n(t,p)\leq \inf\limits_{h>0} \e^{-hx}\E \e^{h(\eta -\l)}=
\exp\Bgl\{ x -(x+\l) \ln \ffrac{x+\l}\l\Bgr\} ,
 \tag 1.18$$
where $x= tn$, $\l = n\s^2$, and $\eta$ is a Poisson random variable with parameter $\l$.
It is  shown in  the proof of Theorem 1.1 in B ~2004,  that if
$T_n$ is a sum of $n$ independent copies of a Bernoulli random variable $\ve=\ve(\s^2)$,
then
$$ \inf_{h<x} (t-h)^{-2} \E (T_n-h)_+^{2}\leq
 \inf_{h<x} (t-h)^{-2} \E (\eta -\l -h)_+^{2},\tag 1.19$$
 where $\eta$ is a Poisson random variable  with parameter $\l =n\s^2$. The right hand side of (1.19)
 is given as an explicit function of $\l$ and $x$ in BKZ~ 2006.

\bigskip
{\bf Remark } {\bf 1.7}. A law of transformation $\{\s_1^2,\dots ,\s_n^2\}\mapsto \s^2$ in (1.1), (1.3) and (1.4)
is a linear function. In bounds involving skewness and kurtosis corresponding transformations
are non-linear, see, for example, (1.12), where the transformation
$\{\ga_1,\dots ,\ga_n\}\mapsto \ga$
is given explicitly.

\bigskip

\bigskip

\bigskip

\bigskip

\head 2. An analytic
description of   transforms $G_2$ of  binomial survival functions $G$.
 \endhead

\bigskip
In this section we recall an explicit analytical description of the right hand side of (1.4)
$$G_2(x) \ \=\ \inf_{h<x} (x-h)^{-2} \E (T_n-h)_+^2,$$
where $T_n$ is a sum of $n$ independent copies of the Bernoulli random variable (1.2).
The description is taken from BKZ ~2006.
Let
 $G(x)=\P\left\{T_{n}\geq x\right\}$ be the survival function
of  $T_{n}$.
 The probabilities $p$, $q$ and the variance $\s^2$ are defined in ~(1.2).  Write
$\lambda=pn$.
The sum $T_{n}=\varepsilon_{1}+\cdots+\varepsilon_{n}$ assumes the values
$$d_{s}=-n\sigma^2+s(1+\sigma^2)\equiv \ffrac{s-\lambda}{q},\quad s=0,1,...,n.$$
The related probabilities satisfy
$$p_{n,s}=\P\left\{T_{n}=d_{s}\right\}=\bbinom{n}{s}q^{n-s}p^{s}.$$
The values $G(d_{s})$ of the survival function of the random variable $T_{n}$ are given by
$$G(d_{s})=p_{n,s}+\cdots+p_{n,n}.$$
Write
$$\nu_{n,s}=\ffrac{sp_{n,s}}{G(d_s)}.$$
Now we can  describe the transform $G_{2}$. Consider a sequence
$$0=r_{0}<r_{1}< \ldots <r_{n-1}<r_{n}=n$$
of points which divide the interval $[0,n]$ into $n$ subintervals $[r_{s},r_{s+1}]$.
To define $G_{2}$ take
$$r_s=\ffrac{\lambda-p\nu_{n,s}}{q\nu_{n,s}+\lambda-s}, \text{\quad $s=0,1, \ldots ,n-1$,}$$
and
$$G_2(x)=\ffrac{\lambda+\nu_{n,s}(s-\lambda-p)-q\nu^2_{n,s}}{qx^2-2q\nu_{n,s}+
\lambda+\nu_{n,s}(s-\lambda-p)}G(d_s), \quad r_s\leq x\leq r_{s+1}.$$

\bigskip

\bigskip

\bigskip

\head 3. Proofs \endhead

\bigskip

{\it Proof of Theorem} {\bf 1.2}. This theorem is a special case of Theorem 1.4.
Indeed, choosing
$$ \s_k^2=\infty ,\quad
 \vk_k=\infty ,\quad k=1,2,\dots , n, $$
we have $v(\vk_k) =\infty$.  Hence  $\alpha_k^2$ from the condition of Theorem 1.4 have to satisfy
$\alpha_k^2 \geq   u^2(\ga_k)$. We choose $\alpha_k^2 =   u^2(\ga_k)$. Then $\s^2=
 ( u^2(\ga_1) +\dots +u^2(\ga_n) )/n$. A small calculation shows that with such $\s^2$
 the skewness $\ga=\s -\ffrac 1 {\s} $ of Bernoulli random variables $\ve =\ve(\s^2)$
 coincides with the expression given in
(1.12) ~$\square$

\bigskip

{\it Proof of Theorem} {\bf 1.3}. This theorem is a special case of Theorem 1.4.
Indeed, choosing
$$ \s_k^2=\infty ,\quad
 \ga_k=- \infty ,\quad k=1,2,\dots , n, $$
we have $u(\ga_k) =\infty$.  Hence  $\alpha_k^2$ from the condition of Theorem 1.4 have to satisfy
$\alpha_k^2 \geq   v(\vk_k)$. We choose $\alpha_k^2 =   v(\vk_k)$. Then $\s^2=
 ( v(\vk_1) +\dots +v(\vk_n) )/n$.  ~$\square$

\bigskip
In the proof of Theorem 1.4 we use the next two lemmas.

\bigskip

\proclaim{Lemma 3.1} Assume that a random variable $X\leq 1$ has mean $\E X=0$, variance $s^2 =\E X^2$,
and skewness such that $\ffrac {\E X^3}{s^3}\geq g$. Then
$$ s^2\leq u^2(g),\quad \quad u(x)\  \= \ \sqrt{1+\ffrac {x^2}4}-\ffrac x2.
\tag 3.1$$
\endproclaim

\demo{Proof} It is clear that
$$(t+s^2)^2(1-t)\geq 0 \quad \text{for} \ t\leq 1 .\tag 3.2 $$
Replacing in (3.2) the variable $t$ by $X$ and taking the expectation, we get $s^2-s^4\geq \E X^3$.
Dividing by $s^3$ and using $\ffrac {\E X^3}{s^3}\geq g$, we derive $\ffrac 1s-s\geq g$.
Elementary considerations show that the latter inequality implies (3.1).
~$\square$\enddemo

\bigskip

\bigskip

\proclaim{Lemma 3.2} Assume that a random variable $X\leq 1$ has mean $\E X=0$, variance $s^2 =\E X^2$,
and kurtosis  such that $\ffrac {\E X^4}{s^4}\leq c$ with some $c\geq 1$.
 Then
$$ s^2\leq v(c),\quad \quad   2v(t) \ \=\ t+1+ \sqrt {(t +1)^2-4}.\tag 3.3$$
\endproclaim

\demo{Proof}
 By H\"older's inequality we have $\ffrac {\E X^4}{s^4}\geq 1$. Hence, the condition $c\geq 1$ is natural.
The function $v$ satisfies $v(c)\geq 1$ for $c\geq 1$. Therefore  in cases where $s^2\leq 1$, inequality (3.3)
turns to the trivial $s^2\leq 1\leq v(c)$. Excluding this trivial case from the further considerations,
we assume that $s^2>1$. Write $a=2\s^2-1$.
 Then $a\geq 1$.
It is clear that
$$(t+s^2)^2(1-t)(a-t)\geq 0 \quad \text{for} \ t\leq 1 .\tag 3.4 $$
Replacing in (3.4) the variable $t$ by $X$ and taking the expectation, we get
 $ \E X^4\geq s^2-s^4+s^6$.
Dividing by $s^4$ and using $\ffrac {\E X^4}{s^4}\leq c$, we derive $\ffrac 1{s^2}-1+s^2\leq c$.
Elementary considerations show that the latter inequality implies (3.3).
~$\square$\enddemo

\bigskip

{\it Proof of Theorem} {\bf 1.4}.
The proof starts with an application of the Chebyshev inequality similar to (1.7).
This reduces the estimation of $\P \{M_n\geq x\}$ to estimation of expectations
$$\E \exp \{ h M_n\}, \quad \E   (M_n-h)_+^2 .$$
As it is noted in the proof of Lemma 4.4 in B~ 2004, it suffices to estimate
$\E (M_n-h)_+^2 $
since the desired bound for the other expectation is implied by
$$\E (M_n-h)_+^2\leq \E (T_n-h)_+^2.\tag 3.5$$
Let us prove (3.5). By Lemma 3.1 the condition $g_k\geq \ga_k$ implies $s^2\leq u^2(\ga_k)$.
While applying Lemma 3.1 one has to replace $X$ by $X_k$, etc. In a similar way, by Lemma 3.2
the condition ~${c_k\leq \vk_k}$ implies $s_k^2\leq v(\vk_k)$. Combining the inequalities and the assumption $s_k^2\leq \s_k^2$,
 we have
$$s_k^2\leq  \min\{ \s_k^2,\, u^2(\ga_k),\, v(\vk_k)\}.\tag 3.6$$
The inequality (3.6) together with the condition of the theorem yields $s_k^2\leq \alpha_k^2$.
As it is shown in the proof of Theorem 1.1 in B~ 2004, the latter inequality
implies (3.5).
~$\square$



\enddocument
The following section of this message contains a file attachment
prepared for transmission using the Internet MIME message format.
If you are using Pegasus Mail, or any another MIME-compliant system,
you should be able to save it or view it from within your mailer.
If you cannot, please ask your system administrator for assistance.

  ---- File information -----------
    File:  unboundedmethod-2008jan24.dvi
    Date:  27 Jan 2008, 12:07
    Size:  105952 bytes.
    Type:  Unknown

The following section of this message contains a file attachment
prepared for transmission using the Internet MIME message format.
If you are using Pegasus Mail, or any another MIME-compliant system,
you should be able to save it or view it from within your mailer.
If you cannot, please ask your system administrator for assistance.

  ---- File information -----------
    File:  unboundedmethod-2008jan24.tex
    Date:  27 Jan 2008, 12:07
    Size:  81672 bytes.
    Type:  Unknown

3 attachments — Download all attachments  
bentkus-juskevicius-skewness-2008-ams.dvi	bentkus-juskevicius-skewness-2008-ams.dvi
48K   Download  
unboundedmethod-2008jan24.dvi	unboundedmethod-2008jan24.dvi
104K   Download  
unboundedmethod-2008jan24.tex	unboundedmethod-2008jan24.tex
80K   Download  
Reply
		
Forward

Reply
More|
[Offline] Tomas Juskevicius	
The following section of this message contains a file attachment prepared for...
	
1/28/08
[Offline] Tomas JuskeviciusLoading...	
1/28/08
[Offline] Tomas Juskevicius to Paulius
	
show details 1/28/08
	
- Show quoted text -

---------- Forwarded message ----------
From: Bentkus@ktl.mii.lt <Bentkus@ktl.mii.lt>
Date: Jan 27, 2008 1:54 PM
Subject: perskaityti!
To: Tomas Juskevicius <tomas.juskevicius@gmail.com>

Aciu uz laiskus!

Gal ir geriau kad busi Lietuvoje, galima bus uzsiimti mokslu.

Man atrodo kad gal geriau eiti pas Rackauska, gal jis
maziau trukdys mokytis.

As pridedu du straipsnius.
Vienas musu, as perrasiau, galvoju atiduoti i LMJ.
Galima buti atspausdinti ir kitur, bet LMJ garantuotai iseis
uz 3 men. Tau reikia perskaityti teksta,
surasti visas klaidas. Tex faile
siulomus pataisymus inesk naudodamas bold arba bf, salia
jau esamo teksto, kad as lengvai matyciau.  Sita reikia
padaryti greitai, geriau per 1 ar 2 dienas.

Kitas straipsnis yra mano, as ji irgi atiduodu i LMJ.
Ten yra nemazai klaideliu. Ryt atsiusiu naujesni
varianta. Ten yra atsakymai i kai kurious tave dominancius klausimus!

Laukiu reakcijos!
Best, V

The following section of this message contains a file attachment
prepared for transmission using the Internet MIME message format.
If you are using Pegasus Mail, or any another MIME-compliant system,
you should be able to save it or view it from within your mailer.
If you cannot, please ask your system administrator for assistance.

  ---- File information -----------
    File:  bentkus-juskevicius-skewness-2008-ams.dvi
    Date:  27 Jan 2008, 12:00
    Size:  48528 bytes.
    Type:  Unknown

\comment Bentkus and Juskevicius skewness Jan 2008
\endcomment

\input amstex
\input epsf

\TagsOnRight
\magnification=\magstephalf

\def\E{{\,{\Bbb E}\,}}

\def\P{{\,{\Bbb P}\,}}
\def\di={\overset{\text{${\Cal D }$}}\to =}
\define\llapp#1{\hbox{\llap{${#1}$}}}

\define\rlapp#1{\hbox{\rlap{${#1}$}}}

\define\na{\,\, {\raise.4pt\hbox{$\shortmid$}}{\hskip-2.0pt\to}\, \, }

\define\<{\bigl <}\define\8{\left<}\define\9{\right>}
\define\>{\bigr>}
\define\ovln#1{\,{\overline{\!#1}}}

\define\PM{^{\2\prime}}

\redefine\D #1{{\<D{#1},{#1}\>}}
\define\norm#1{\left\|#1\right\|}
\define\nnorm#1{\bgl\|#1\bgr\|}
\define\snorm#1{\bigl\|\2#1\2\bigr\|}
\define\ssqrt{^{1/2}}


\define\weht#1{\bigl(\2\e+\nnorm{#1}\t\2\bigr)}
\define\leei#1#2{\log\E e^{i\8\right.\!{#1},{#2}\!\left.\9}}
\define\eei#1#2{\E e^{i\8\right.\!{#1},{#2}\!\left.\9}}
\define\lee#1#2{\log\E e^{\8\right.\!{#1},{#2}\!\left.\9}}
\define\ee#1#2{\E e^{\8\right.\!{#1},{#2}\!\left.\9}}
\define\xxxi#1{{\overline\xi}_{#1}(h)}

\define\oxh#1#2{{\overline\xi_{#1}}(h_{#2})}
\define\ooxh#1#2{{\overline\xi{}_{#1}\PM}(h'_{#2})}
\define\xxi#1{\xi_#1}
\define\4{\hskip1pt}

\define\wthj#1{\widetilde h_{#1}}

\define\twhtj#1{{\th\4\nnorm{\wthj#1}\tau}}

\define\xwthj#1{{\overline\xi_{{#1}}}(\wthj{#1})}
\define\xhpj#1#2{{\overline\xi{}\PM_{{#1}}}(h'_{#2})}

\define\detDj#1{(2\4\pi)^{d/2}\,(\det D_{#1})\ssqrt}

\define\2{\hskip.5pt}
\define\bgr#1{\4\bigr#1}
\define\bgl#1{\bigl#1\4}
\define\3{\hskip-10pt}
\define\R#1#2{\Cal R_d\bgl(#1,#2\bgr)}
\define\RB#1#2{\Cal R_d\Bgl(#1,#2\Bgr)}
\define\RBM#1#2{\Cal R_{d-1}\Bgl(#1,#2\Bgr)}

\define\pfrac#1#2{{}^{#1}\hskip-.5pt\!\big/\!\hskip.5pt{}_{#2}}
\define\sfrac#1#2{{#1}\big/{#2}}
\define\ffrac#1#2{\raise.5pt\hbox{\eightpoint$\4\dfrac{\,#1\,}{\,#2\,}\4$}}
\define\Bgr#1{{\eightpoint\4\raise.5pt\hbox{{$\Bigr#1$}}}}
\define\Bgl#1{{\eightpoint\raise.5pt\hbox{{$\Bigl#1$}}\4}}
\define\step#1{\raise9pt\hbox{\eightpoint$\tsize{#1}$}}

\documentstyle {amsppt}
\hoffset  0.25  true cm
\voffset  1.5 true cm

\hsize=16 true cm
\vsize=22 true cm
\loadmsbm
\def\t{\theta}

\def\D{\Delta }
\def\ve{\varepsilon }

\def\vk{\varkappa}
\def\ga{\gamma }

\def\l{\lambda }

\def\s{\sigma }

\def\s{\sigma}
\def\t{\theta}

\def\={\overset{ \text{def} }\to =}

\def\R{\Bbb R}

\def\Z{\Bbb Z}

\def\I{\Bbb I}

\define\e{\operatorname{e}}

\def\fs#1#2#3{ {#1}_{#2},\dots ,  {#1}_{#3}    }
\def\fsu#1#2#3{ {#1}_{#2}+\dots + {#1}_{#3}    }

\define\bbinom#1#2{\raise.5pt\hbox{\eightpoint$\dsize\4
\binom{\raise-1.3pt\hbox{${\ssize\,#1\,}$}}
{\raise1.3pt\hbox{${\ssize \,#2\,}$}}\4$}}

\documentstyle{amsppt}

\topmatter
\title Bounds for tail probabilities of martingales
using skewness and kurtosis
\endtitle
\rightheadtext{ Skewness and curtosis  } \leftheadtext{V. Bentkus}
\author V. Bentkus\footnotemark"$^{1}$" and  T. {Ju\v skevi\v cius}$^{1}$
\endauthor

\footnotetext"$^{1}$"{ Vilnius Institute of Mathematics and
Informatics. Supported by grant ???}

\address
Vilnius institute of mathematics and informatics,\ 
Akademijos 4,\ 
LT-08663  Vilnius
\endaddress

\comment
and \smallskip
Vilnius pedagogical university \smallskip
Studentu 39\smallskip
LT-08106  Vilnius \smallskip
Lithuania  \smallskip
\endcomment

\email
bentkus\@ktl.mii.lt, ?????????????
\endemail

\date January 2008
\enddate
\keywords
skewness, kurtosis,
Hoeffding's
inequalities,
sums of independent random variables, martingales,
bounds for tail probabilities, Bernoulli random variables, binomial
tails, probabilities of large deviations, method of bounded
differences
  \endkeywords
\subjclass
60E15
\endsubjclass

\abstract
Let $M_n= \fsu X1n$ be a sum of independent random variables such that
$ X_k\leq 1$,  $\E X_k =0$ and $\E X_k^2=\s_k^2$ for all $k$.
Hoeffding 1963, Theorem 3, proved that
$$\P\{M_n \geq nt \}\leq H^n(t,p ) ,\quad H(t,p)= \bgl( 1+qt/p\bgr)^{p +qt}
\bgl( {1-t}\bgr)^{q -qt}$$
with
$$q=\ffrac 1{1+\s^2},\quad p=1-q, \quad \s^2 =\ffrac { \s_1^2+\dots
+\s_n^2}n,\quad 0<t<1.$$
Bentkus 2004 improved Hoeffding's
inequalities using binomial tails as upper bounds.
Let $\ga_k =\E X_k^3/\s_k^3$ and $ \vk_k= \E X_k^4/\s_k^4$ stand for the skewness
and kurtosis of $X_k$.
In this paper we prove (improved) counterparts of the Hoeffding inequality
replacing $\s^2$ by certain functions of $\fs \ga 1n$ respectively $\fs \vk 1n$.
Our bounds extend to a general setting where $X_k$ are martingale
differences, and they
can combine
the knowledge of skewness and/or kurtosis and/or variances of ~$X_k$.
Up to  factors
bounded by  $e^2/2$ the bounds are final.
All our results are new since no inequalities incorporating
skewness or kurtosis control so far are known.
\endabstract

\endtopmatter

\document


\head  1. Introduction and results
\endhead

In a celebrated paper of Hoeffding 1963 several inequalities for sums of bounded random
variables were established. For improvements of the Hoeffding inequalities and
related results see, for example,   Talagrand ~1995, McDiarmid 1989,
Godbole and Hitczenko~ 1998,  Pinelis ~1998--2007,
Laib ~1999, B ~2001--2007,
van de Geer ~2002, Perron~ 2003,
BGZ ~2006--2006, BGPZ ~2006,
BKZ ~2006, ~2007, BZ 2003, etc.
Up to certain constant factors, these improvements are close to the final
optimal inequalities, see ~B ~2004, ~BKZ ~2006.
However so far no  bounds taking into account  information related to
skewness and/or kurtosis are known, not to mention
certain  results related to symmetric random variables, see  BGZ 2006, BGPZ 2006.
In this paper we prove general and optimal counterparts of
Hoeffding's 1963 Theorem 3, using assumptions related to skewness and/or  kurtosis.

\bigskip

Let us recall Hoeffding's 1963 Theorem 3.
Let $M_n =\fsu X1n$ be a sum of independent random variables such that
$ X_k\leq 1$, $\E X_k=0$, and $\E X_k^2=\s_k^2$ for all $k$.
Write
$$ \s^2=\ffrac {\s_1^2+\dots +\s_n^2}n, \quad p=\ffrac {\s^2} {1+\s^2},\quad q=1-p .$$
Hoeffding ~1963, Theorem 3, established the inequality
$$\P\{M_n \geq nt \}\leq H^n(t,p ) ,\quad H(t,p)= \bgl( 1+qt/p\bgr)^{p +qt}
\bgl( {1-t}\bgr)^{q -qt}  \tag 1.1$$
assuming that $ 0<t <1$.
One can rewrite $ H^n(x,p ) $ as
$$ H^n(x,p )=\inf_{h>0} \exp\{ -h nt \}\E \exp\{hT_n\},$$
where $T_n=\fsu \ve 1n $ is a sum of $n$
independent copies of a  Bernoulli random variable, say ~${\ve=\ve(\s^2)}$,
such that
$$\P\{\ve =-\s^2 \}=q,\quad   \P\{\ve=1\}=p,\quad \E \ve^2 =\s^2.\tag 1.2 $$
Using the shorthand $x=nt$, we can rewrite the Hoeffding result as
$$\P\{M_n \geq x \}\leq \inf_{h>0} \e^{ -h x }\E \e^{hT_n}, \tag 1.3$$
In B 2004 the inequality (1.3) is improved to
$$\P\{M_n \geq x \}\leq \inf_{h<x} (x-h)^{-2} \E (T_n-h)_+^2
.\tag 1.4 $$
Actually, inequalities  (1.1), (1.3) and (1.4) extend to  cases where $M_n$ is a martingale or even super-martingale,
see B~ 2004 for a proof.
In  the case of (1.1) and (1.3) this was noted already by Hoeffding 1963.

The right hand side of (1.4) satisfies
$$ \inf_{h<x} (t-h)^{-2} \E (T_n-h)_+^{2}\leq \ffrac {\e^2}2\P \{ T_n\geq x \}, \quad \e=2.718... \tag 1.5 $$
for integer $x\in \Z$. For non-integer $x$ one has to interpolate the probability
$\log$-linearly, see ~B ~2004 for details. The right-hand side of (1.4) can be given explicitly
as a function of $x$, $p$ and $n$,
see~ BKZ 2006, as well as Section 2 of the present paper. To have
bounds as tight as possible is essential for statistical applications, like those in audit,
see BZ 2003.

Our intention in this paper
is to develop methods leading to counterparts of (1.1), (1.3) and ~(1.4) such that
information related to the skewness and kurtosis
$$\ga_k= \ffrac {\E ( X_k-\E X_k)^3}{\s_k^3},\quad \quad \vk_k= \ffrac {\E (X_k-\E X_k)^4}{\s_k^4}\tag 1.6$$
of $X_k$ is taken into account
(in this paper we define  $\ga_k=\infty$ and $\vk_k=1$
if $\s_k=0$).
  All our results hold in general martingale setting.

All known proofs of inequalities of type (1.3) and (1.4) start with an application of
Chebyshev's inequality. For example, in the case of (1.4)
we can estimate
$$\P\{M_n \geq x \}\leq \inf_{h<x} (x-h)^{-2} \E (M_n-h)_+^2
\tag 1.7 $$
since the indicator function $t\mapsto \I \{ t\geq x\}$ obviously satisfies
$\I \{ t\geq x\}\leq (x-h)^{-2}  (t-h)_+^2$ for all ~${t\in \R}$.
The further proof of (1.4) consists in showing that
$\E (M_n-h)_+^2\leq \E (T_n-h)_+^2$ for all ~${h\in \R}$. We would like to emphasize that
all our proofs are optimal in the sense that no further improvements are possible
in estimation of $\E (M_n-h)_+^2$. Indeed, in the special case $M_n=T_n$
the inequality $\E (M_n-h)_+^2\leq \E (T_n-h)_+^2$ turns into the equality
$\E (T_n-h)_+^2= \E (T_n-h)_+^2$.

In view of (1.7) it is  natural to introduce and to study
transforms $G\mapsto G_\beta$
of survival functions $G(x)=\P \{X\geq x\}$
of the type
$$G_\beta (x) = \inf_{h<x} (x-h)^{-\beta}\E (X-h)_+^\beta ,  \quad  \beta>0,\tag 1.8$$
defining $G_0=G$ in the case $\beta =0$. See Pinelis ~1988, 1989, B~ 2004, ~BKZ 2006 for related known results.

The paper is organized as follows. In the Introduction we provide necessary definitions and formulations
of  our results, including their versions for sums of martingale differences.
In Section 2 we recall a description of the transform $G\mapsto G_2$ of binomial survival
functions---our bounds are given using $G_2$.
Section 3 contains proofs of the results.

\bigskip
Henceforth $M_n=\fsu X1n$ stands for  a martingale sequence such that the differences ~$X_k$
are uniformly  bounded (we set $M_0=X_0=0$).
Without loss of generality we can assume that the bounding constant is $1$,
that is, that  $X_k\leq 1$. Let $\Cal F_0 \subset \Cal F_1\subset \dots \subset \Cal F_n$
be a related sequence of $\sigma$-algebras
such that $M_k$ are $\Cal F_k$-measurable.
Introduce the conditional variance $s_k^2$, skewness ~$g_k$ and kurtosis $c_k$ of $X_k$ by
$$s_k^2= \E \bgl(X_k^2\, \Bgl|\, \Cal F_{k-1}\bgr),\quad
g_k= \E \bgl(  {X_k^3}\, \Bgl|\, \Cal F_{k-1}\bgr)/{s_k^3} ,\quad
c_k= \E \bgl( {X_k^4}\, \Bgl|\, \Cal F_{k-1}\bgr)/{s_k^4}.\tag 1.9
$$
Note that $s_k^2$, $g_k$, $c_k$ are $\Cal F_{k-1}$-measurable random variables.

\bigskip

{\bf Remark} {\bf 1.1}. We prove our results using (1.4) for martingales. It is proved in B ~2004 that
all three inequalities ~(1.1), ~(1.3) and
~(1.4) hold with $\s^2=(\s_1^2+\dots +\s_n^2)/n$
if $M_n$ is a martingale with differences ~${X_k\leq 1 }$ such that the conditional variances $s_k^2$
satisfy $s_k^2\leq \s_k^2$ for all $k$.

\bigskip

It is easy to check that Bernoulli random variables $\ve=\ve(\s^2)$ of type (1.2) have  variance $\s^2$ and
skewness ~$\ga$
related as
$$\ga = \s -\ffrac 1\s ,\quad \s^2 = u^2(\ga ),\quad \text{where} \ u(x) =\sqrt{1+\ffrac {x^2}4}-\ffrac x2
.\tag 1.10 $$

\bigskip

\proclaim{Theorem 1.2} Assume that the differences $X_k$ of a martingale $M_n$
satisfy ~${X_k\leq 1}$, and that the conditional skewness $g_k$
of $X_k$ are bounded from below by some non-random $\ga_k$, that is, that
$$ g_k \geq \ga_k,\quad k=1,2,\dots , n. \tag 1.11$$
Then $(1.3)$ and $(1.4)$ hold
with $T_n$ being  a sum of $n$ independent copies of a Bernoulli random variable ~${\ve=\ve(\s^2)}$
of type $(1.2)$ with
skewness $\ga$ and variance $\s^2$ defined by
$$\ga= \ffrac 1{\sqrt{n}}\ffrac{\ga_1 u(\ga_1) +\dots + \ga_n u (\ga_n) }{\sqrt{u^2(\ga_1)+ \dots +u^2(\ga_n)}},
\quad \s^2 = \ffrac{ u^2(\ga_1)+\dots +u^2(\ga_n) }n. \tag 1.12 $$
In the special case where all $\ga_k$ are equal, $\ga_1=\dots =\ga_n=\ga$, the Bernoulli random variable has
skewness $\ga$ and variance $\s^2= u^2(\ga)$.
\endproclaim

\bigskip

It is easy to see that Bernoulli random variables $\ve=\ve(\s^2)$ of type (1.2) have  variance $\s^2$ and
kurtosis ~$\vk$
related as
$$
\vk  = \ffrac 1{\s^2}-1+\s^2 ,\quad 2\s^2=  {\vk+1\pm \sqrt {(\vk +1)^2-4}}
.\tag 1.13 $$
In particular
$$
\s^2\leq v(\vk ),
\quad \text{where} \ 2v(t) =t+1+ \sqrt {(t +1)^2-4}.\tag 1.14$$

\bigskip
\proclaim{Theorem 1.3}
Assume that the differences $X_k$ of a martingale $M_n$
satisfy ~${X_k\leq 1}$, and that the conditional kurtosis $c_k$
of $X_k$ are bounded from above by some non-random $\vk_k$, that is, that
$$ c_k \leq \vk_k,\quad k=1,2,\dots , n. \tag 1.15$$
Then $(1.3)$ and $(1.4)$ hold
with $T_n$ being  a sum of $n$ independent copies of a Bernoulli random variable ~${\ve=\ve(\s^2)}$
of type $(1.2)$ with
kurtosis  $\vk$ and variance $\s^2$ defined by
$$\vk  = \ffrac 1{\s^2}-1+\s^2,
\quad \s^2 = \ffrac{ v(\vk_1)+\dots +v(\vk_n) }n,\tag 1.16 $$
where the function $v$ is given in $(1.14)$.
In the special case where  $\vk_1=\dots =\vk_n=\vk$, the Bernoulli random variable has
kurtosis  $\vk$ and variance $\s^2= v(\vk)$.
\endproclaim

\bigskip
The next Theorem 1.4 allows to combine our knowledge about variances, skewness and kurtosis. Theorems 1.2, 1.3
and (1.4), (1.4) for martingales (see Remark 1.1)
are special cases of Theorem 1.4 setting in various combinations $\s_k^2=\infty$, $\ga_k=-\infty$, $\vk_k=\infty$.

\bigskip
\proclaim{Theorem 1.4}
Assume that the differences $X_k$ of a martingale $M_n$
satisfy ~${X_k\leq 1}$, and that their
conditional variances $s_k^2$, skewness $g_k $ and
kurtosis $c_k$ satisfy
$$s_k^2\leq \s_k^2,\quad g_k \geq \ga_k,\quad
c_k \leq \vk_k,\quad k=1,2,\dots , n \tag 1.17$$
with some non-random $s_k^2\geq 0$, $g_k\geq -\infty$ and $1\leq c_k\leq \infty$.
Assume that  numbers $\alpha_k^2$ satisfy
$$\alpha_k^2 \geq \min\{ \s_k^2,\,  u^2(\ga_k),\,  v(\vk_k)\}.$$
Then $(1.3)$ and $(1.4)$ hold
with ~${T_n}$ being  a sum of $n$ independent copies of a Bernoulli random variable ~${\ve=\ve(\s^2)}$
of type $(1.2)$ with
$$ \s^2 = \ffrac{ \alpha_1^2 +\dots +\alpha_n^2 }n,$$
where  functions $u$ and  $v$ are defined  in $(1.10)$ and $(1.14)$ respectively.
\endproclaim

\bigskip
{\bf Remark } {\bf 1.5}. All our  inequalities can be extended to the case where $M_n$ is a super-martingale.
Furthermore, their maximal versions hold, that is, in the left hand sides of these inequalities we
can replace $\P\{M_n\geq x\}$ by $\P\left\{\max\limits_{1\leq k \leq n}M_k\geq x\right\}$.

\bigskip
{\bf Remark } {\bf 1.6}. One can estimate the right hand sides of our inequalities
using Poisson distributions. In the case of Hoeffding's functions this is done by Hoeffding ~1963.
In notation of (1.1) his bound is
$$ H^n(t,p)\leq \inf\limits_{h>0} \e^{-hx}\E \e^{h(\eta -\l)}=
\exp\Bgl\{ x -(x+\l) \ln \ffrac{x+\l}\l\Bgr\} ,
\tag 1.18$$
where $x= tn$, $\l = n\s^2$, and $\eta$ is a Poisson random variable with parameter $\l$.
It is  shown in  the proof of Theorem 1.1 in B ~2004,  that if
$T_n$ is a sum of $n$ independent copies of a Bernoulli random variable $\ve=\ve(\s^2)$,
then
$$ \inf_{h<x} (t-h)^{-2} \E (T_n-h)_+^{2}\leq
\inf_{h<x} (t-h)^{-2} \E (\eta -\l -h)_+^{2},\tag 1.19$$
where $\eta$ is a Poisson random variable  with parameter $\l =n\s^2$. The right hand side of (1.19)
is given as an explicit function of $\l$ and $x$ in BKZ~ 2006.

\bigskip
{\bf Remark } {\bf 1.7}. A law of transformation $\{\s_1^2,\dots ,\s_n^2\}\mapsto \s^2$ in (1.1), (1.3) and (1.4)
is a linear function. In bounds involving skewness and kurtosis corresponding transformations
are non-linear, see, for example, (1.12), where the transformation
$\{\ga_1,\dots ,\ga_n\}\mapsto \ga$
is given explicitly.

\bigskip

\bigskip

\bigskip

\bigskip

\head 2. An analytic
description of   transforms $G_2$ of  binomial survival functions $G$.
\endhead

\bigskip
In this section we recall an explicit analytical description of the right hand side of (1.4)
$$G_2(x) \ \=\ \inf_{h<x} (x-h)^{-2} \E (T_n-h)_+^2,$$
where $T_n$ is a sum of $n$ independent copies of the Bernoulli random variable (1.2).
The description is taken from BKZ ~2006.
Let
$G(x)=\P\left\{T_{n}\geq x\right\}$ be the survival function
of  $T_{n}$.
The probabilities $p$, $q$ and the variance $\s^2$ are defined in ~(1.2).  Write
$\lambda=pn$.
The sum $T_{n}=\varepsilon_{1}+\cdots+\varepsilon_{n}$ assumes the values
$$d_{s}=-n\sigma^2+s(1+\sigma^2)\equiv \ffrac{s-\lambda}{q},\quad s=0,1,...,n.$$
The related probabilities satisfy
$$p_{n,s}=\P\left\{T_{n}=d_{s}\right\}=\bbinom{n}{s}q^{n-s}p^{s}.$$
The values $G(d_{s})$ of the survival function of the random variable $T_{n}$ are given by
$$G(d_{s})=p_{n,s}+\cdots+p_{n,n}.$$
Write
$$\nu_{n,s}=\ffrac{sp_{n,s}}{G(d_s)}.$$
Now we can  describe the transform $G_{2}$. Consider a sequence
$$0=r_{0}<r_{1}< \ldots <r_{n-1}<r_{n}=n$$
of points which divide the interval $[0,n]$ into $n$ subintervals $[r_{s},r_{s+1}]$.
To define $G_{2}$ take
$$r_s=\ffrac{\lambda-p\nu_{n,s}}{q\nu_{n,s}+\lambda-s}, \text{\quad $s=0,1, \ldots ,n-1$,}$$
and
$$G_2(x)=\ffrac{\lambda+\nu_{n,s}(s-\lambda-p)-q\nu^2_{n,s}}{qx^2-2q\nu_{n,s}+
\lambda+\nu_{n,s}(s-\lambda-p)}G(d_s), \quad r_s\leq x\leq r_{s+1}.$$

\bigskip

\bigskip

\bigskip

\head 3. Proofs \endhead

\bigskip

{\it Proof of Theorem} {\bf 1.2}. This theorem is a special case of Theorem 1.4.
Indeed, choosing
$$ \s_k^2=\infty ,\quad
\vk_k=\infty ,\quad k=1,2,\dots , n, $$
we have $v(\vk_k) =\infty$.  Hence  $\alpha_k^2$ from the condition of Theorem 1.4 have to satisfy
$\alpha_k^2 \geq   u^2(\ga_k)$. We choose $\alpha_k^2 =   u^2(\ga_k)$. Then $\s^2=
( u^2(\ga_1) +\dots +u^2(\ga_n) )/n$. A small calculation shows that with such $\s^2$
the skewness $\ga=\s -\ffrac 1 {\s} $ of Bernoulli random variables $\ve =\ve(\s^2)$
coincides with the expression given in
(1.12) ~$\square$

\bigskip

{\it Proof of Theorem} {\bf 1.3}. This theorem is a special case of Theorem 1.4.
Indeed, choosing
$$ \s_k^2=\infty ,\quad
\ga_k=- \infty ,\quad k=1,2,\dots , n, $$
we have $u(\ga_k) =\infty$.  Hence  $\alpha_k^2$ from the condition of Theorem 1.4 have to satisfy
$\alpha_k^2 \geq   v(\vk_k)$. We choose $\alpha_k^2 =   v(\vk_k)$. Then $\s^2=
( v(\vk_1) +\dots +v(\vk_n) )/n$.  ~$\square$

\bigskip
In the proof of Theorem 1.4 we use the next two lemmas.

\bigskip

\proclaim{Lemma 3.1} Assume that a random variable $X\leq 1$ has mean $\E X=0$, variance $s^2 =\E X^2$,
and skewness such that $\ffrac {\E X^3}{s^3}\geq g$. Then
$$ s^2\leq u^2(g),\quad \quad u(x)\  \= \ \sqrt{1+\ffrac {x^2}4}-\ffrac x2.
\tag 3.1$$
\endproclaim

\demo{Proof} It is clear that
$$(t+s^2)^2(1-t)\geq 0 \quad \text{for} \ t\leq 1 .\tag 3.2 $$
Replacing in (3.2) the variable $t$ by $X$ and taking the expectation, we get $s^2-s^4\geq \E X^3$.
Dividing by $s^3$ and using $\ffrac {\E X^3}{s^3}\geq g$, we derive $\ffrac 1s-s\geq g$.
Elementary considerations show that the latter inequality implies (3.1).
~$\square$\enddemo

\bigskip

\bigskip

\proclaim{Lemma 3.2} Assume that a random variable $X\leq 1$ has mean $\E X=0$, variance $s^2 =\E X^2$,
and kurtosis  such that $\ffrac {\E X^4}{s^4}\leq c$ with some $c\geq 1$.
Then
$$ s^2\leq v(c),\quad \quad   2v(t) \ \=\ t+1+ \sqrt {(t +1)^2-4}.\tag 3.3$$
\endproclaim

\demo{Proof}
By H\"older's inequality we have $\ffrac {\E X^4}{s^4}\geq 1$. Hence, the condition $c\geq 1$ is natural.
The function $v$ satisfies $v(c)\geq 1$ for $c\geq 1$. Therefore  in cases where $s^2\leq 1$, inequality (3.3)
turns to the trivial $s^2\leq 1\leq v(c)$. Excluding this trivial case from the further considerations,
we assume that $s^2>1$. Write $a=2\s^2-1$.
Then $a\geq 1$.
It is clear that
$$(t+s^2)^2(1-t)(a-t)\geq 0 \quad \text{for} \ t\leq 1 .\tag 3.4 $$
Replacing in (3.4) the variable $t$ by $X$ and taking the expectation, we get
$ \E X^4\geq s^2-s^4+s^6$.
Dividing by $s^4$ and using $\ffrac {\E X^4}{s^4}\leq c$, we derive $\ffrac 1{s^2}-1+s^2\leq c$.
Elementary considerations show that the latter inequality implies (3.3).
~$\square$\enddemo

\bigskip

{\it Proof of Theorem} {\bf 1.4}.
The proof starts with an application of the Chebyshev inequality similar to (1.7).
This reduces the estimation of $\P \{M_n\geq x\}$ to estimation of expectations
$$\E \exp \{ h M_n\}, \quad \E   (M_n-h)_+^2 .$$
As it is noted in the proof of Lemma 4.4 in B~ 2004, it suffices to estimate
$\E (M_n-h)_+^2 $
since the desired bound for the other expectation is implied by
$$\E (M_n-h)_+^2\leq \E (T_n-h)_+^2.\tag 3.5$$
Let us prove (3.5). By Lemma 3.1 the condition $g_k\geq \ga_k$ implies $s^2\leq u^2(\ga_k)$.
While applying Lemma 3.1 one has to replace $X$ by $X_k$, etc. In a similar way, by Lemma 3.2
the condition ~${c_k\leq \vk_k}$ implies $s_k^2\leq v(\vk_k)$. Combining the inequalities and the assumption $s_k^2\leq \s_k^2$,
we have
$$s_k^2\leq  \min\{ \s_k^2,\, u^2(\ga_k),\, v(\vk_k)\}.\tag 3.6$$
The inequality (3.6) together with the condition of the theorem yields $s_k^2\leq \alpha_k^2$.
As it is shown in the proof of Theorem 1.1 in B~ 2004, the latter inequality
implies (3.5).
~$\square$



\enddocument